\documentclass [12pt]{article}
\usepackage {amssymb}
\begin {document}
\sloppy

\title
{Stationary quantum stochastic processes from the cohomological
point~of~view}

\author
{Grigori G. Amosov\footnote {The work was supported by
INTAS-00-738}\\
Department of Higher Mathematics\\
Moscow Institute of Physics and Technology\\
Dolgoprudni 141700\\
RUSSIA\\
E-mail: gramos@deom.chph.ras.ru}

\maketitle

\abstract
{Stationary quantum stochastic process $j$ is introduced as a
*-homomorphism embedding an involutive graded algebra
$\tilde K=\oplus _{i=1}^{\infty }K_i$ into a ring of
(abelian) cohomologies of the
one-parameter group $\alpha $ consisting of *-automorphisms of certain
operator algebra in a Hilbert space such that every $x$ from $K_i$ is
translated into an additive $i-\alpha$-cocycle $j(x)$.
It is shown that (noncommutative) multiplicative
markovian cocycle defines a perturbation of the stationary
quantum stochastic process in the sense of such definition. The
$E_0$-semigroup $\tilde \beta $ on the von Neumann algebra $\cal N$
associated with the markovian perturbation of $K$-flow $j$ posseses
the restriction $\tilde \beta |_{{\cal N}_0},\ {\cal N}_0\subset
{\cal N},$
which is conjugate to the flow of Powers shifts $\beta $ associated
with $j$. It yields for $\tilde \beta $ an analogue of the Wold
decomposition for classical stochastic process
on completely nondeterministic and deterministic parts.
The examples of quantum stationary stochastic processes on
the algebras of canonical commutation,
anticommutation and square of white noise relations are considered.
In the model situation of the space $L^2(\mathbb R)$
all markovian cocycles of the group of shifts are described
up to unitary
equivalence of perturbations.}

\section {Introduction.}

Let $\nu $ be a measure on the real axe $\mathbb R$.
Denote $S=(S_t)_{t\in \mathbb R}$ the flow of shifts
acting on the measurable functions $f$ on $\mathbb R$ by
the formula $(S_tf)(x)=f(x-t),\ x,t\in \mathbb R$.
Then given a number $r\in \mathbb R$
the function $I_r(t)=\nu ([r,r+t]),\ t\in \mathbb R$,
is a $1-S$-cocycle, i.e. $I_r(t+s)=I_r(t)+S_t(I_r(s)),\ t,s\in
\mathbb R$.  Analogously fixing numbers $r_1,\dots ,r_k\in \mathbb R$
we get that the function $I_{r_1\dots r_k}(t_1,\dots
,t_k)= I_{r_1}(t_1)S_{t_1}(I_{r_2}(t_2)S_{t_2}(\dots
I_{r_k}(t_k)\dots ),\ t_i\in \mathbb R$, associated with the
tensor product $\nu ^{\otimes k}$ is a $k-S$-cocycle with
the characteristic property $S_{t_1}(I(t_2,\dots ,t_{k+1}))-
I(t_1+t_2,t_3,\dots ,t_k)+\dots +(-1)^iI(t_1,\dots ,t_{i-1},
t_i+t_{i+1},t_{i+2},\dots ,t_{k+1})+\dots +(-1)^kI(t_1,\dots ,t_k)=0,\
t_i\in \mathbb R$. One can consider a ring of cohomologies
$H^*=\oplus H^i$ generated by the measure $\nu $ in this way.
The canonical bilinear map
$(x_i,x_j)\to x_i\cup x_j\in H^{i+j},\ x_i\in H^i,x_j\in H^j,$
defining a ring structure in $H^*$ can be obtained from the
action of $S$ by the formula $(x_i\cup x_j)(t_1,\dots ,t_{i+j})
=x_i(t_1,\dots ,t_i)S_{t_1+\dots +t_i}(x_j(t_{i+1},\dots ,t_{i+j}))$
(see \cite {Gui80}).
We consider a
graded algebra of cohomologies $A=\oplus _{i=1}^{\infty }A_i$ such
that $A_i$ consists of additive $i$-cocycles of the group of
automorphisms $\alpha $ associated with stationary quantum stochastic
process $j$ over an involutive algebra $K$. In our construction $A_1$
is generated by the basic operator-valued stochastic measures which
are the creation, anihilation and number of particles processes in
the important applications. Given an involutive algebra $K$ we put
$K_1=K$ and construct the graded algebra $\tilde K= \oplus
_{i=1}^{\infty }K_i$ with respect to a certain linear
associative operation $\odot :K_i\times K_j\to K_{i+j}$. In the case
when $K$ is a Lie algebra one can choose the universal enveloping
algebra for $\tilde K$.  We define a stationary
quantum stochastic process as a *-homomorphism $j$ of $\tilde K$ into
the standard ring of (abelian) cohomologies
of the group $\alpha $ such that every $x\in K_i$ is maped to
an $i-\alpha$-cocycle $j(x)$.  The homomorphism $j$ transferes
the operation $\odot $ in $\tilde K$ to the
cohomological multiplication $\cup$.  We
consider a class of cocycle perturbations of $j$ by markovian
cocycles. The markovian cocycle perturbation of $K$-flow
constructed through our procedure determines
the associated $E_0$-semigroup $\tilde \beta$ on the von Neumann
algebra $\cal N$ which posseses the restriction $\tilde \beta |_{{\cal
N}_0},\ {\cal N}_0\subset {\cal N},$ conjugate to the flow
of Powers shifts associated with the initial $K$-flow.
Using the technics of \cite {Bul2}, it is possible to
extract an automorphic part of the $E_0$-semigroup.
Hence our result defines an analogue of the Wold
decomposition for classical stochastic process on completely
nondeterministic and deterministic parts. We give several examples
where $j$ determines quantum stochastic processes on the algebra of
canonical commutation or anticommutation relations and the square of
white noise correspondingly.  The basic ideas are illustrated on the
model of markovian perturbations for the group of shifts in
$L^2(\mathbb R)$.  In this case the associated markovian cocycles are
constructed in the explicit form by means of the inner function
techniques.

\section {Markovian perturbations of the group of shifts in
$L^2(\mathbb R)$.}

Let $S=(S_t)_{t\in \mathbb R}$ be a strong continuous group of
unitary operators on a Hilbert space $\cal H$.
A strong continuous family of unitaries
$W=(W_t)_{t\in \mathbb R}$ in $\cal H$ is called $a\ multiplicative\
1-S-cocycle$ if $W_{t+s}=W_tS_tW_sS_{-t},\ s,t\in \mathbb R,\ W_0=I$.
The cocycle $W$ is said to be $a\ multiplicative\
1-S-coboundary$ if there exists a unitary opertor $J$ defining
$W$ by the formula $W_t=JS_tJ^*S_{-t},\ t\in \mathbb R$.
Every multiplicative cocycle $W$ determines
a new unitary group $U=(U_t)_{t\in \mathbb R}$ in $\cal H$ by the
formula $U_t=W_tS_t,\ t\in \mathbb R$. This group can be named a
cocycle perturbation of $S$.
Notice that if $W$ is a coboundary, then $W_tS_t=JS_tJ^*,\
t\in \mathbb R,$ i.e. the coboundary determines the perturbation
which is unitary equivalent to the initial group.
Consider the group of shifts $S=(S_t)_{t\in \mathbb R}$ acting in
the Hilbert space
${\cal H}=L^2(\mathbb R)$ by the formula $(S_tf)(x)=f(x+t),\ f\in
{\cal H}$.
The group of cohomologies $H^1(S,L^2(\mathbb R))$ is generated by
$an\ additive\ 1-S-cocycle$ $\chi =(\chi _t)_{t\in \mathbb R}$
defined by $\chi _t(x)=1,\ -t<x\leq 0,\ \chi _t(x)=0$ otherwise. The
cocycle $\chi$ satisfies the characteristic properties
$\chi_{t+s}=\chi_t+S_t\chi_s,\ ||\chi_t-\chi_s||=2|t-s|^{1/2},\
t,s\in \mathbb R$. Let ${\cal H}_t$ be a subspace of $\cal
H$ generated by all functions with supports belonging to the segment
$[-t,+\infty )$.  Notice that ${\cal H}_t=S_t{\cal H}_0,\ t\in
\mathbb R$.  We shall call the multiplicative cocycle $W$ by
$markovian$ if $W_t|_{{\cal H}\ominus {\cal H}_t}=I,\ t\geq 0$. This
property means that $W$ doesn't perturbe "a future" of the system.
The Markov property for perturbations was introduced in \cite {Acc}.
Using the cocycle property
$W_{-t}=S_tW_t^*S_{-t},\ t\in \mathbb R,$ we can rewrite
the condition of markovianity in the form
$W_{-t}|_{{\cal H}\ominus {\cal H}_0}=I,\ t\geq 0$.
It garantees that linear operators
$V_t=U_{-t}|_{{\cal H}_0},\ t\geq 0,$ are correctly defined and
form a $C_0$-semigroup of isometries $V$ in the Hilbert space
${\cal H}_0$. We call $V$ by a semigroup
associated with the markovian cocyle perturbation of $S$.
Given a $C_0$-semigroup of isometrical operators $V$
in the Hilbert space ${\cal H}_0$, one can define the Wold
decomposition ${\cal H}_0={\cal H}^{(0)}\oplus {\cal H}^{(1)}$
on the subspace ${\cal H}^{(0)}$ reducing $V$ to the semigroup of
unitary operators and the subspace ${\cal H}^{(1)}$ reducing
$V$ to the semigroup of completely nonunitary operators isomorphic
to the semiflow of right shifts in ${\cal K}\otimes L^2({\mathbb
R_+})$, where $\cal K$ is a Hilbert space with the dimension
equal to the deficiency index of the
generator of $V$ (see \cite {Nik}).  The semigroups $V|_{{\cal
H}^{(0)}}$ and $V|_{{\cal H}^{(1)}}$ can be named a unitary part
and a shift part of the semigroup $V$ correspondingly.

{\bf Proposition 2.1.} {\it The deficiency index of the generator of
the semigroup $V$ associated with the markovian cocycle
perturbation of $S$ equals $1$.}

{\bf Proof.}

Consider a family of functions $\zeta_t=W_{-t}\chi_{-t},\ t\geq 0$,
where $\chi$ is the additive $1-S$-cocycle defined in the begin of
the section. It follows from the definitions of $\chi$ and $W$ that
the family $\zeta$ is continuous.  Then
$\zeta_{t+s}=W_{-t-s}\chi_{-t-s}=W_{-t}S_{-t}W_{-s}S_t(\chi_{-t}+S_{-t}
\chi_{-s}) =
W_{-t}\chi_{-t}+W_{-t}S_{-t}W_{-s}\chi_{-s}=\zeta_t+V_t\zeta_s,\
t,s\geq 0$, where we have used the identity
$W_{-s}S_t\chi_{-t}=S_t\chi_{-t},\ s,t\geq 0$, following from
the markovian property of $W$ in the form
$W_{-s}|_{{\cal H}\ominus {\cal H}_0}=I$.
Hence $\zeta=(\zeta_t)_{t\geq 0}$
is an additive $1-V$-cocycle.  Notice that
$(V_{tn}\zeta_t,\zeta_t)=(W_{-t}S_{-t}V_{t(n-1)}\zeta_t,
W_{-t}\chi_{-t})= (S_{-t}V_{t(n-1)}\zeta_t,\chi_{-t})=0,\
n\in \mathbb N ,\ t>0$.  Therefore the sum $\xi _t=\sum \limits
_{n=0}^{+\infty } e^{-tn}V_{tn}\zeta_t$ is well defined and we can
write the integral sum $\xi =\lim \limits _{t\to 0}\xi _t=\int
\limits _{0}^{+\infty }e^{-t}d\zeta_t$. It follows that $V_t^*\xi
=e^{-t}\xi $. So we have proved that the deficiency index of the
generator of $V$ more or equal to one. Let ${\cal H}_{\zeta}$ be
a subspace of ${\cal H}$ generated by $U_s\zeta _t,\ s,t\geq 0$.
The cocycle property $\zeta _{t+s}=\zeta_t+V_t\zeta _s=
\zeta _t+U_{-t}\zeta _s,\ s,t\geq 0,$ leads to
the invariance of ${\cal H}_{\zeta}$ under the action of $U_t,\
t\in \mathbb R$. Moreover $U_t\zeta _t=W_tS_tW_{-t}\chi _{-t}=
W_tS_tW_{-t}S_{-t}S_t\chi _{-t}=S_t\chi _{-t}=-\chi _t,\ t\geq 0$.
Hence ${\cal H}\ominus {\cal H}_{0}\subset
{\cal H}_{\zeta}$. Put ${\cal H}^{(1)}={\cal H}_{\zeta}\cap {\cal
H}_0$, then the subspace ${\cal H}^{(0)}={\cal H}_0\ominus {\cal
H}^{(1)}$ is invariant under the action of $U_t,\ t\in \mathbb R$.
Hence the restriction $V|_{{\cal H}^{(0)}}$ consists of unitary
operators. We get the Wold decomposition ${\cal H}_0={\cal
H}^{(0)}\oplus {\cal H}^{(1)}$ for the semigroup of isometries $V$.
Therefore the index of $V$ equals the index of $V|_{{\cal H}^{(1)}}$
which is one.

{\bf Theorem 2.2.} {\it Let $W$ be a markovian cocycle, then
there  exists $s-\lim \limits _{t\to +\infty }W_{-t}=W_{-\infty }$
and an isometrical operator $W_{-\infty }$ defines the Wold
decomposition ${\cal H}_0={\cal H}^{(0)}\oplus {\cal H}^{(1)}$
for the semigroup $V$ associated with the markovian perturbation
by $W$
such that ${\cal H}^{(1)}=W_{-\infty }{\cal H}_0$.}

{\bf Proof.}

Firstly we shall apply the technics which is anlogues to one in the
proof of theorem 2.5 in \cite {Bhat}.
Notice that $W_{-t-s}f=W_{-s}S_{-s}W_{-t}S_sf=
W_{-s}f$ for all $f\in {\cal H}\ominus {\cal H}_{-s},\
s,t\geq 0$ by the markovian property. Hence the sequence
$W_{-t}f$ converges when $t$ tends to $+\infty $ for the dense
set of vectors $f$. Therefore the limit exists by the Banach-Steinhaus
theorem. One can see that the subspace
${\cal H}^{(1)}=\vee _{t\geq 0}W_{-t}({\cal H}_0\ominus {\cal
H}_{-t})$ is generated by the functions $\zeta _t=W_{-t}\chi _{-t},
\ t\geq 0,$ which form an additive $1-V$-cocycle with the
orthogonal increaments (see the proof of Proposition 2.1).
It follows that the restriction of $V$ to ${\cal H}^{(1)}$
is unitary equivalent to the semiflow of right shifts in
$L^2(\mathbb R)$.

Let $\cal V$ be a subspace of ${\cal H}_0$ invariant under
the action of the semigroup of right shifts $S'=S^*|_{{\cal H}_0}$.
Then (see \cite {Nik}) ${\cal V}=M_{\Theta }H_0$, where ${\cal
F}^{-1}M_{\Theta }{\cal F}$ is the operator of multiplication
by an inner function $\Theta $, $\cal F$ is the Fourier transform.
Notice that $M_{\Theta }$ is an isometrical operator.  Denote
$P_{[0,t]}$ and $P_{[t,+\infty )}$ the orthogonal projections on the
spaces ${\cal H}_0\ominus {\cal H}_{-t}$ and ${\cal H}_{-t},\ t\geq
0,$ correspondingly. Given a unitary $C_0$-group $R=(R_t)_{t\in
\mathbb R}$ in the Hilbert space ${\cal H}_0\ominus {\cal V}$, define
a family of unitary operators in ${\cal H}_0$ by the formula
$$
W_{-t}=(R_tP_{{\cal H}_0\ominus {\cal V}}S_t+P_{\cal V})P_{[t,+\infty
)}+ M_{\Theta } P_{[0,t]},\ t\geq 0.
\eqno (2.1)
$$
In the following proposition we introduce the model describing all
markovian cocycles up to unitary equivalence of perturbations.
To be exact, for every markovian cocycle $\tilde W$
there exists the markovian cocycle $W$ of the form given in
the proposition such that $\tilde W_t=J_tW_t,\ t\in \mathbb R,$
where the $1-U$-coboundary $J_t=JW_tS_tJ^*W_{-t}S_{-t},\ t\in \mathbb
R,$ of the group $U=(W_tS_t)_{t\in \mathbb R}$
is defined by the unitary operator $J$ satisfying the relation
$\tilde W_tS_t=JW_tS_tJ^*,\ t\in \mathbb R$.

{\bf Proposition 2.3.} {\it The family $W$ given by (2.1)
in ${\cal H}_0$ and acting
identicaly in ${\cal H}
\ominus {\cal H}_0$ defines
a multiplicative markovian cocycle such that
$\lim \limits _{t\to +\infty }W_{-t}f =M_{\Theta }f$ for
$f\in {\cal H}_0$ and $\lim \limits _{t\to +\infty }W_{-t}f=f$ for
$f\in {\cal H}\ominus {\cal H}_0$. The unitary part of the semigroup
$V$ associated with the markovian cocycle perturbation by $W$ is
$R$.}

{\bf Proof.}

Extend the family $W$ defined in (2.1) for $t\geq 0$ by the
formula $W_t=S_tW_{-t}^*S_{-t},\ t\geq 0$. Consider
the set of unitary operators $U_t=W_tS_t,\ t\in \mathbb R$.
Notice that $U_{-t}M_{\Theta }=M_{\Theta }S_{-t},\
U_{-t}f\in {\cal V},\ f\in {\cal H}_t,\ t\geq 0$. Hence the subspace
${\cal L}={\cal V}\oplus ({\cal H}\ominus {\cal H}_0)$ is invariant
under action of $U_t,\ t\in \mathbb R,$ and the restriction $U|_{\cal
L}$ is unitary equivalent to the group $S$.  To complete the proof
notice that the restriction $U$ to the subspace ${\cal H}_0\ominus
{\cal V}={\cal H}\ominus {\cal L}$ coincides with $R^*$.

{\bf Theorem 2.4.}
{\it Given a unitary group $\tilde R$ which is uniformly
continuous or has a pure point spectrum,
there exist the inner function $\Theta $ and the unitary group
$R$ in the Hilbert space ${\cal H}_0\ominus M_{\Theta }{\cal H}_0$
unitary equivalent to $\tilde R$ such that the
markovian cocycle $(2.1)$ satisfies the condition
$W_t-I\in s_2,\ t\in \mathbb R$.}

Theorem 2.4 for the case of uniformly continuous $\tilde
R$ can be found in \cite {Amo1},\cite {Amo2}. The proof for
$\tilde R$ with a pure point spectrum is also constructed in
cited papers in the implicit form.
One can compare the condition $W-I\in s_2$ on the unitary operator
$W$ appearing in the theorem
with the Feldman criterion on the equivalence of Gaussian measures
(see \cite {Fel},\cite {Gui72}) and the Araki criterion of
the quasi-equivalence for quasifree states of the algebra of
canonical commutation relations
(see \cite {ArYa}).
It seems that the unitary operators $W$ satisfying
our condition translate equivalent states one to another.
It is
also useful to notice that the condition $W_t-I\in s_2,\ t\in \mathbb
R,$ is nessesary and sufficient for the family $W=(W_t)_{t\in \mathbb
R}$ to define an inner cocycle on the hyperfinite factor generated by
a quasifree representation of the algebra of canonical
anticommutation relations (see \cite {Amo2},\cite {MuYa}).

\section {Stationary quantum stochastic process as *-homomorphism
into a ring of cohomologies.}

Let $K$ be an involutive algebra and $K_i=K^{\otimes i}$.
Supply the family $(K_i)_{i=1}^{+\infty }$
with a linear associative operation
$\odot $ defining left and right actions of every $x\in K_i$
on $K_j$ such that $x\odot y,\ y\odot x\in K_{i+j},\ y\in K_j$.
We assume that $K_i\odot K_j=K_{i+j},\ i,j\in \mathbb N$.
So we obtained the graded algebra $\tilde K=\oplus _{i=1}^{\infty }
K_i$ with respect to the multiplication defined by the operation
$\odot $. If $K$ is a Lie algebra, it is possible to take
the unversal enveloping algebra for $\tilde K$ and the
multiplication in $\tilde K$ for $\odot $.
Consider a one-parameter
$w^*$-continuous group of *-automorphisms
$\alpha =(\alpha _t)_{t\in {\mathbb R}}$ on
the algebra $\cal B(H)$ of all bounded operators in a
Hilbert space $\cal H$.
Suppose that the action of $\alpha $ can be correctly
defined on certain involutive algebra $\cal M$ consisting of
linear operators (in general unbounded) in $\cal H$.
Consider the standard resolvent for $\alpha $ constructed from
nonhomogeneous chains such that (see \cite {Gui80})
$$
0\longrightarrow {\cal M}\longrightarrow
Hom(\mathbb R,{\cal M})\stackrel{d_1}{\longrightarrow }\dots
\longrightarrow Hom(\mathbb R^i,{\cal M})\stackrel{d_i}
{\longrightarrow}\dots
$$
where $d_i(x)(t_1,\dots ,t_{i+1})=\alpha _{t_1}(x(t_2,\dots
,t_{i+1}))-x(t_1+t_2,\dots ,t_{i+1})+\dots +(-1)^ix(t_1,\dots ,t_i),
\ x=x(t_1,\dots ,t_i)\in Hom(\mathbb R^i,{\cal M})$.
Denote $A_i=kerd_i/Imd_{i-1}$, then $A=\oplus _{i=1}^{+\infty }A_i$
is a ring with respect to the multiplication
defined by a bilinear map $(x,y)\stackrel {\cup }\to x(t_1,\dots ,t_i)
\alpha _{t_1+\dots +t_i}(y(t_{i+1},\dots ,t_{i+j}))\in
A_{i+j},\ x\in A_i,\ y\in A_j$.

Every *-homomorphism $j$ of the graded algebra
$\tilde K$ into the graded algebra $A$ such that every $x\in K_i$ is
translated into an $i-\alpha $-cocycle $j(x)$ we shall call
by {\it a stationary quantum stochastic process} over the algebra
$K$. Notice that
our definition is based upon the well-known one given in \cite {Acc1}.
It is also useful to remark that we don't need $0-\alpha $-cohomologies
in our construction. Sometimes
we can recognize two processes $j$ and $\tilde j$
determining the cocycles $j(x)$ and $\tilde j(x)$
differ on the coboundary for the
fixed $x$ as obtaining one from other by a shift in time.
For example, given a stationary quantum stochastic process $j$
the $1-\alpha $-cocycle $j(x)(t)$ is differ on the coboundary
$\alpha _t(j(x)(r))-j(x)(r)$ from
the $1-\alpha $-cocycle $j(x)(r+t)-j(x)(r)$ which is associated
with the stationary quantum stochastic process $\tilde j$
obtained from $j$ by a shift in time on $r$.
In applications we claim that $j$ keeps the
basic algebraic structure in $K$ and the mutiplication $\odot $ in
$\tilde K$ but we do not need to require for $j$ the preserving of the
multiplication in $K$ if it is defined (see the next section).
We shall suppose that the operators involved in the
image of $j$ generate whole $\cal M$.

{\bf Proposition 3.1.} {\it Let $j$ be defined on
$K=K_1$. Then there exists a unique extension of $j$ to
whole $\tilde K$.}

{\bf Proof.}

Notice that $j$ translates the operation $\odot $ in
the cohomological multiplication $\cup $.
Denote $j^{(i)}=j|_{K_i}$, then
one can obtain the action of $j$ by the induction,
$$
j^{(k+l)}_{t_1,\dots ,t_{k+l}}(x\odot y)=
j^{(k)}_{t_1,\dots ,t_k}(x)\alpha _{t_1+\dots +t_k}
(j^{(l)}_{t_{k+1},\dots ,t_{k+l}}(y)).
$$
for all $x\in K_k,\ y\in K_l$. Here we use the property
$K_i\odot K_j=K_{i+j}$.

\section {Stationary quantum stochastic processes on the algebras of
canonical commutation relations, the square of white noise relations
and canonical anticommutation relations.}

Consider an involutive Lie algebra $K$ generated by the elements
$B,\ B^+,\ \Lambda ,\ {\bf 1}$ satisfying the relations
$[B,B^+]={\bf 1},\ [\Lambda ,B]=-B,\ [\Lambda ,B^+]=B^+$.
Let the operation $\odot $ is generated by the multiplication
in the universal enveloping algebra of $K$. One can define a
stationary quantum stochastic process $j$ over $K$ by the formula
$$
j(B)=b_t,\ j(B^+)=b_t^*,
$$
$$
j(\Lambda )=\Lambda _t,\ j({\bf 1})=t{\bf 1},
$$
where $b_t,b_t^*,\Lambda _t$ are the bosonic anihilation, creation
and number of particles processes (see \cite {HuPa}).

Analogously if $K$ is the Lie algebra generated by
the elements $B^-,B^+,M$ satisfying the relations of $\bf sl_2$
which are $[B^-,B^+]=M,\ [M,B^{\pm}]=\pm 2B^{\pm}$, one can
consider the universal enveloping algebra of $K$ and we obtain
the quantum Levy process generated by the square of white noise (SWN)
constructed in \cite {Acc2}, $$ j(B^-)=b_t,\ j(B^+)=b^+_t, $$ $$
j(M)=\gamma t+n_t,
$$
where the basic processes $b_t,b_t^+$ and $n_t$ satisfy the relations
of SWN,
$$
b_tb_t^+-b_t^+b_t=\gamma t+n_t,\
n_tb_t-b_tn_t=-2b_t,
$$
$$
n_tb_t^+-b_t^+n_t=2b_t^+,\
(b_t)^*=b_t^+,\ n_t^*=n_t,
$$
with a fixed parameter $\gamma >0$ and $t\in
\mathbb R$.

The algebra of canonical anticommutation relations (CAR) is
generated by elements $a_t,a^*_t,\ t\in \mathbb R,$ satisfying
the relations $a_ta^*_s+a^*_sa_t=t\wedge s{\bf 1},\
a_ta_s+a_sa_t=0$.
A graded algebra $K=\oplus _{i=1}^{\infty }K_i$
can be obtained in the following
way, $K_i$ is a tensor product of
i-th copies of the algebra of $2\times 2$-matrix units.
Let elements $a_k,a_k^*,\ k\in \mathbb N,$ satisfy the
canonical anticommutation relations $a_ka_l^*+a_l^*a_k=\delta
_{kl}{\bf 1},\ a_ka_l+a_la_k=0,\ k,l\in \mathbb N$.
Then $K_i$ is generated by $a_k,a_k^*,\ 1\leq k\leq i$.
Determine a canonical
operation $\odot $ by the formula $x_1\odot x_2\odot \dots \odot x_n=
y_1y_2\dots y_n$, where $y_{i}=a_{i},\
a^*_{i},\ a^*_{i}a_{i},\ a_ia^*_i$ if
$x_i=a,\ a^*,\ a^*a,\ aa^*$ correspondingly.
Then a stationary quantum stochastic process $j$ can be
defined by the formula
$$
j(a)=a_t,\ j(a^*)=a^*_t,
$$
$$
j(a^*a)=\Lambda _t,\ j({\bf 1})=t{\bf 1},
$$
where $a_t,a_t^*,\Lambda _t$ are the basic Fermion processes
(see \cite {ApHu}). Notice that $j$ satisfies the relations
$j(a^*)=j(a)^*,\ j(a^*a+aa^*)=j(a)^*j(a)+j(a)j(a)^*=j({\bf 1}),\
j([a^*a,a])=[j(a^*a),j(a)]=-j(a)=a_t,\
j([a^*a,a^*])=[j(a^*a),j(a^*)]=j(a^*)=a_t^*$, but it is not the
algebraic *-morphism because $j(a^*a)\neq j(a)^*j(a)$.

\section {Cocycle perturbations of $K$-flows and the Wold
decomposition.}

We shall use the notation of previous parts of this paper.
Remember that strong continuous family $W=(W_t)_{t\in {\mathbb R}}$
of unitary operators in $\cal H$ is named a multiplicative $\alpha
$-cocycle if $W_{t+s}=W_t\alpha _t(W_s),\ s,t\in \mathbb R$.  Suppose
that the action of $W$ is correctly defined on $\cal M$, i.e.
$W_txW_t^*$ are well defined for all $t\in \mathbb R,\ x\in {\cal
M}$. Let ${\cal M}_{t]}$, ${\cal M}_{[s}$ and ${\cal M}_{[s,t]}$ be
involutive subalgebras of $\cal M$ generated by all increaments of
the form $j(x)(r)-j(x)(l),\ x\in K_1,$ where $l\leq r\leq t$, $s\leq
l\leq r$ and $s\leq l\leq r\leq t$ correspondingly. We shall call
(see also \cite {Acc}) a multiplicative $\alpha $-cocycle by
{\it markovian} (with respect to
the stationary quantum stochastic process $j$) if $W_t{\cal
M}_{t]}W_t^*\subset {\cal M}_{t]}$ and $W_txW_t^*=x$ for all $x\in
{\cal M}_{[t},\ t\geq 0$.

{\bf Theorem 5.1.} {\it For any markovian cocycle $W$
the formula
$$
\tilde j(x)(t)=j(x)(t),\ t\geq 0,
$$
$$
\tilde j(x)(t)=W_tj(x)(t)W_t^*,\ t\leq 0,
$$
$$
\tilde \alpha _t(\cdot )=
W_t\alpha _t(\cdot )W_t^*,\ x\in K_1,\ t\in \mathbb R,
$$
defines a new stationary quantum
stochastic process $\tilde j$ over $K$ with an associated group
of automorphisms $\tilde \alpha $.}

{\bf Proof.}

Analogously to the proof of Proposition 2.1 we
obtain
$$
\tilde j(x)(t+s)=W_{t+s}j(x)(t+s)W_{t+s}^*=
W_{t+s}(j(x)(t)+\alpha _t(j(x)(s)))W_{t+s}^*=
$$
$$
W_t\alpha _t(W_s)j(x)(t)\alpha _t(W_s^*)W_t^*+
W_t\alpha _t(W_sj(x)(s)W_s^*)W_t^*=
$$
$$
W_tj(x)(t)W_t^*+\tilde \alpha _t(W_sj(x)(s)W_s^*)=
\tilde j(x)(t)+\tilde \alpha _t(\tilde j(x)(s)),
$$
$s,t\leq 0$.
Here we used the identity
$\alpha _t(W_s)j(x)(t)\alpha _t(W_s^*)=
j(x)(t)$ due to the markovian property
$W_s\alpha _{-t}(j(x)(t))W_s^*=-\alpha _s(W_{-s}^*)j(x)(-t)
\alpha _s(W_{-s})=-\alpha _s(W_{-s}^*(j(x)(-t-s)-j(x)(-s))W_{-s})=
-\alpha _s(j(x)(-t-s)-j(x)(-s))=-j(x)(-t)=\alpha _{-t}(j(x)(t)),\
s,t\leq 0$. One can extend $\tilde j(x)(t)$ for $t\geq 0$
using the cocycle condition for $\tilde j(x)(t)$. It yealds
$\tilde j(x)(t)=j(x)(t),\ t\geq 0$.
To complete the proof we only need
to apply Proposition 3.1.

{\bf Proposition 5.2.} {\it
$$
\tilde \alpha _t(x)=\alpha _t(x),\ x\in {\cal M}_{[0},\ t\geq 0,
$$
$$
\tilde \alpha _{-t}(x)=\alpha _{-t}(x),\ x\in {\cal M}_{[t},\ t\geq
0.
$$
}

{\bf Proof.}

It immidiately follows from the markovian property of $W$ that
$\tilde \alpha _t(x)=W_t\alpha _t(x)W_t^*=\alpha _t(x),\
x\in {\cal M}_{[0},\ t\geq 0$.
The markovian property implies that $W_txW_t^*=x,\ x\in {\cal
M}_{[t},\ t\geq 0,$ which is
equivalent to $\alpha _t(W_{-t})^*x\alpha _t(W_{-t})=x,\
x\in {\cal M}_{[t}, t\geq 0,$ or $W_{-t}^*xW_{-t}=x,\
x\in {\cal M}_{[0},\ t\geq 0$,
by the cocycle condition for $W$. Hence
$\tilde \alpha _{-t}(x)=W_{-t}\alpha _{-t}(x)W_{-t}^*=
\alpha _{-t}(x),\ x\in {\cal M}_{[t},\ t\geq 0$.

Denote ${\cal N=M}''\cap {\cal B(H)},\
{\cal N}_{t]}={\cal M}_{t]}''\cap {\cal B(H)},\
{\cal N}_{[t}={\cal M}_{[t}''\cap {\cal B(H)}$ and ${\cal
N}_{[s,t]}={\cal M}_{[s,t]}''\cap {\cal B(H)}$ the corresponding von
Neumann algebras. Notice that ${\cal N}_{t+s]}=\alpha _t({\cal
N}_{s]}),\ t,s\in \mathbb R$.  We shall call a stationary quantum
stochastic process $j$ by a K-flow and the group $\alpha $ associated
with $j$ by a group of automorphisms associated with K-flow if the
following conditions hold, $$ {\cal N}_{s]}\subset {\cal N}_{t]},\
t>s, $$ $$ \vee _{t\in {\mathbb R}}{\cal N}_{t]}={\cal N}, $$ $$
\wedge _{t\in
{\mathbb R}}{\cal N}_{t]}={\bf C1}
$$
(see \cite {Emc}).

{\bf Proposition 5.3.} {\it
If there exists a
vector $\Omega \in \cal H$ which is cyclic and separating with
respect to $\cal N$ and the increaments of a stationary
quantum stochastic process $j$
are independent in the classical (commutative) sence that $\phi
(xy)=\phi (x)\phi (y),\ x\in {\cal N}_{[t},\ y\in {\cal N}_{t]},\
t\in {\mathbb R},$ for the state $\phi (\cdot )=(\Omega ,\cdot \Omega
)$, then $j$ is a $K$-flow.  }

{\bf Proof.}

Choose $x\in\wedge _{t\in {\mathbb R}}{\cal N}_{t]}$, then
$\phi ((x-\phi (x){\bf 1})y)=\phi (x-\phi (x){\bf 1})
\phi (y)=0$ for all $y\in \vee _{t\in \mathbb R}{\cal N}_{[t}=
\vee _{t\in \mathbb R}{\cal N}_{t]}=
\cal N$. Hence $(x-\phi (x){\bf 1})\Omega =0$ and
$x=\phi (x){\bf 1}$ as $\Omega $ is cyclic and separating.
The result follows from.

Let von Neumann algebras $\tilde {\cal N}_{t]}$,
$\tilde {\cal N}_{[s,t]}$ and $\tilde
{\cal N}_{[t}$ be associated with the perturbed process $\tilde j$ in
the same way as the algebras ${\cal N}_{t]}$,
${\cal N}_{[s,t]}$ and ${\cal N}_{[t}$ are
associated with the process $j$.

{\bf Proposition 5.4.} {\it Let $j$ and $W$ be a K-flow and
a markovian cocycle correspondingly. Then the markovian
perturbation $\tilde j$ is also K-flow.}

{\bf Proof.}

One can see that the von Neumann algebras generated by the
increaments of $\tilde j$ are $\tilde {\cal N}_{t]}=W_t
{\cal N}_{t]}W_t^*\subset {\cal N}_{t]}$ by the markovian
property. Hence
$\wedge _{t\in {\mathbb R}}\tilde {\cal N}_{t]}={\bf C1}.$
The conditions
$$
\tilde {\cal N}_{s]}\subset \tilde {\cal N}_{t]},\ t>s,
$$
$$
\vee _{t\in {\mathbb R}}\tilde {\cal N}_{t]}=\tilde {\cal N}
$$
are satisfied by the definition.

For a stationary quantum stochastic process $j$ one can name
the $E_0$-semigroup $\beta _t=\alpha _{-t}|_{{\cal N}_{0]}},\ t\geq
0,$ on the von Neumann algebra ${\cal N}_{0]}$ by associated with
$j$.  In the case when $j$ is a $K$-flow, the semigroup $\beta
=(\beta _t)_{t\geq 0}$ is a flow of Powers shifts \cite {Pow}, i.e.
$\wedge _{n\in \mathbb N}\beta _{tn}({\cal N}_{0]})={\bf C1},\ t>0$
(see \cite {Bul1}). Fix a stationary quantum stochastic process $j$
with the associated $E_0$-semigroup $\beta $.  For a cocycle $W$
being markovian with respect to $j$ we shall call the $E_0$-semigroup
$\tilde \beta _t(\cdot )=W_{-t}\beta _t(\cdot ) W_{-t}^*,\ t\geq 0,$
on ${\cal N}_{0]}$ by associated with the markovian perturbation
of the initial process by $W$.  Two $E_0$-semigroups $\beta $ and
$\tilde \beta $ on the von Neumann algebras $\cal A$ and $\tilde
{\cal A}$ correspondingly are called to be conjugate if there exist
two injective *-homomorphisms $\theta :{\cal A}\to \tilde {\cal A}$
and $\theta ^+:\tilde {\cal A}\to {\cal A}$ such that $\beta
_t(x)=\theta ^+\beta _t\theta (x),\ \theta ^+\theta (x)=x,\ \theta
\theta ^+(y)=y,\ x\in {\cal A},\ y\in \tilde {\cal A},\ t\geq 0$.

{\bf Theorem 5.5.} {\it Given a markovian perturbation of
the $K$-flow $j$ with the associated flow of Powers shifts
$\beta $ on the von Neumann algebra ${\cal N}_{0]}$ acting in
the Hilbert space $\cal H$ with a ciclyc vector  $\Omega $, there
exists the von Neumann algebra $\tilde {\cal N}_{0]}\subset
{\cal N}_{0]}$ such that the restriction $\tilde \beta |_{\tilde
{\cal N}_{0]}}$
of the $E_0$-semigroup $\tilde \beta $ associated with the
markovian perturbation is conjugate to $\beta $.}

{\bf Proof.}

Put ${\cal H}_t=[{\cal N}_{[-t}\Omega ],\ t\in \mathbb R$,
then the set ${\cal H}_t,\ t\geq 0,$ is dense in
the Hilbert space ${\cal H}$.
Arguing similarily to the proof of theorem 2.2 one can
obtain that there exists
$s-\lim \limits _{t\to +\infty }W_{-t}=
W_{-\infty}$.
Then the injective *-endomorphism $\theta :{\cal N}_{0]}\to
\tilde {\cal N}_{0]}\subset {\cal N}_{0]}$ given by the formula
$\theta (x)=W_{-t}xW_{-t}^*,\ x\in {\cal N}_{[-t,0]},\
t\geq 0$, is well defined because
$W_{-t-s}xW_{-t-s}^*=W_{-t}\alpha _{-t}(W_{-s})x\alpha _{-t}
(W_{-s}^*)W_{-t}^*=W_{-t}xW_{-t}^*$ for all $x\in {\cal N}_{[-t,0]}$
by the markovian property $W_{-s}yW_{-s}^*=y,\ y\in {\cal N}_{[0},$
which implies that $W_{-s}\alpha _t(x)W_{-s}^*=
\alpha _t(x)$ for all $x\in {\cal N}_{[-t,0]}$.
It follows that $\lim \limits _{t\to +\infty }
W_{-t}xW_{-t}^*f=\lim \limits _{t\to +\infty }W_{-t}xW_{-s}^*f=
W_{-\infty }xW_{-s}^*f$ for all $f\in {\cal
H}_{-s},\ x\in {\cal N}_{0]},\ s\geq 0$. Hence the sequence
$W_{-t}xW_{-t}^*f$ converges when $t$ tends to $+\infty $
for all $f\in {\cal H}_{0]}$ by
the Banach-Steinhaus theorem.
Analogously it is possible
to define the injective *-homomorphism $\theta ^+:\tilde {\cal
N}_{0]}\to {\cal N}_{0]}$ by the formula $\theta
^+(x)=W_{-t}^*xW_{-t},\ x\in \tilde {\cal N}_{[-t,0]},\ t\geq 0$.
Notice that $\theta ^+(x)=W_{-\infty }^*xW_{-\infty },\ x\in \tilde
{\cal N}_{0]}$.  One can see that $\tilde \beta _t(x)=\theta \beta _t
\theta ^+ (x),\ x\in \tilde {\cal N}_{0]}$. This proves the theorem.

Earlier it was investigated the existence of "an automorphic part"
in the quantum dynamical semigroup
which is completely compatible with the faithful
normal state (see \cite {Bul2}).
Theorem 5.5 allows to obtain "a shift part" of
the $E_0$-semigroup obtained by a markovian cocycle perturbation from
the flow of Powers shifts.
So it can be considered as some analogue of the picking out
a completely nondetermenistic part in
the Wold decomposition for the classical stochastic
processes.

\section*{Acknowledgments} {The author is grateful to Professor
Luigi Accardi for kind hospitality during his visit at
Centro Vito Volterra Universita di Roma Tor Vergata where
a part of this work was done.}

\end {document}